%% file: rsdPrePrint.tex
\documentclass[11pt]{article}
\usepackage{amsmath}
\usepackage{amsfonts}
\usepackage{amssymb}
\usepackage{color}
\usepackage{bbm, subfigure, appendix}

\usepackage{tikz}
\usetikzlibrary{arrows}
\usetikzlibrary{shapes}
\usetikzlibrary{backgrounds}
\usetikzlibrary{patterns}
\usepackage{graphicx}
\usepackage{pseudocode}

\numberwithin{equation}{section}

\title{Recursive Schur Decomposition{\footnote{Notice: This manuscript has been authored by UT-Battelle, LLC, under
 Contract No. DE-AC05-00OR22725 with the U.S. Department of Energy. The United States Government retains and the 
 publisher, by accepting the article for publication, acknowledges that the United States Government retains a 
 non-exclusive, paid-up, irrevocable, world-wide license to publish or reproduce the published form of this 
 manuscript, or allow others to do so, for United States Government purposes.}}}
 
\author{Rahul S. Sampath{\footnote{Corresponding author: Rahul S. Sampath (\tt{sampathrs@ornl.gov}).}}, Bobby Philip, Srikanth Allu and Srdjan Simunovic \\
Computer Science and Mathematics Division\\ Oak Ridge National Laboratory\\ One Bethel Valley Road, Oak Ridge, TN 37831}

\begin{document}

\maketitle 

\input{abstract}

\input{introduction}

\input{setup}

\input{algorithm}

\input{results}

\input{conclusions}

\input{acknowledgements} 

{
\bibliographystyle{plain}
\bibliography{refs}
}

\end{document}

%% file: abstract.tex
\begin{abstract}
 In this article, we present a parallel recursive algorithm based on multi-level domain decomposition
  that can be used as a precondtioner to a Krylov subspace method to solve sparse linear
 systems of equations arising from the discretization of partial differential 
 equations (PDEs). We tested the effectiveness of the algorithm on several PDEs using different number of sub-domains 
 (ranging from 8 to 32768) and various problem sizes (ranging from about 2000 to
  over a billion degrees of freedom). We report the results from these tests; the results show 
  that the algorithm scales very well with the number of sub-domains.
\end{abstract}

%% file: introduction.tex
\section{Introduction}
\label{sec:intro} 
 Domain decomposition methods are a class of numerical techniques that are used to solve large scale boundary
 value problems on complex domains using several processors. These algorithms are based on the classic
 principle of ``divide and conquer''. In these methods, the domain of interest is divided into a number of 
 sub-domains and the global problem is reduced to a smaller problem defined only on the interfaces 
 between the sub-domains. The solution to the interface problem provides the necessary and sufficient conditions 
 required to split the global problem into a set of independent problems defined on each sub-domain. These 
 sub-problems are easier (with respect to mesh generation) and cheaper (with respect to storage
 and operation count) to solve compared to the global problem; moreover, they can be solved concurrently. For a 
 good review on different domain decomposition methods, we refer the reader to \cite{toselli2005,
  smith2004, mathew2008} and references therein.
 
 The technique of domain decomposition can be applied to both linear and nonlinear problems but, we only focus on linear problems 
 in this article. In particular, we restrict our attention to sparse linear systems of equations that are
 generated by discretization and perhaps linearization of partial differential equations (PDEs). Furthermore, we are only concerned 
 with discretizations using the finite element method (FEM) although, the methods discussed here could
 potentially be applied to other types of discretizations as well.
 
 The interface problem forms the core of any domain decomposition method. It is typically solved (possibly 
 only approximately) using an iterative algorithm because using a direct method would be too expensive. In fact, 
 the matrix (required for a direct method) corresponding to the interface problem is
 seldom constructed explicitly as it would involve solving the sub-domain problems several times. A major
 challenge faced by domain decomposition methods is the fact that the size of the interface problem grows as
 the number of sub-domains increases. For single-level domain decomposition methods, an increase in the interface
 problem size translates to an increase both in the number of iterations required to solve the problem as well 
 as the cost per iteration; this critically limits the scalability of these methods with respect to the number
 of sub-domains. Two-level domain decomposition methods \cite{dohrmann2007, farhat2001, mandel1993} attempt to 
 address this issue by solving an additional coarse global problem, known as Coarse Grid Correction, coupling all
 the sub-domains at each iteration. This makes the number of iterations (equivalently, the convergence rate of the algorithm) required to solve the 
 problem independent of the number of sub-domains but, the cost of applying the Coarse Grid Correction grows 
 with the number of sub-domains. Besides, the Coarse Grid Correction is problem dependent and it is not 
 straightforward to apply these methods to solve arbitrary PDEs; these methods have mostly only been applied to
 solve structural mechanics problems. Multi-level domain decomposition techniques \cite{henon2006, gaidamour2008} take a different approach - they
 construct a hierarchy of interface problems using hierarchical domain decomposition such that the interface problem 
 size at any level of the hierarchy is independent of the total number of sub-domains. These interface problems 
 are nested, i.e., to solve the interface problem at any level of the hierarachy we need to solve the interface 
 problems at the level immediately below it. It would be too expensive to solve these interface problems exactly; 
 hence, an incomplete factorization is used instead. Due to these approximations, the performance of these methods depends
 on the number of levels in the hierarchy. As the number of subdomains increases so does the number of levels
 in the hierarchy; this typically results in a deterioration in the performance of these algorithms.

 In this work, we use a different technique to address the problem described above. Our algorithm uses a hierarchical 
 domain decomposition framework as well but, it is distinct from existing multi-level methods in many respects. Presently,
 we have only applied this method to solve two dimensional PDEs on long and thin domains. This type of computational domain 
 configuration is common in many problems of significant practical importance. For example, it is used to model
 oil pipelines \cite{joseph1997, paidoussis1998}, electric power lines \cite{zhang2006}, nuclear fuel pins \cite{olander1976, olander2009} 
 and arteries \cite{grinberg2008, manguoglu2008}.
 
 The next section describes how we construct our hierarchical domain decomposition. Our multi-level solver is
 presented in Section \ref{sec:algorithm}. In Section \ref{sec:results}, we report on some numerical experiments
 performed using this algorithm. The paper ends with some concluding remarks.

%% file: setup.tex
\section{Hierarchical Domain Decomposition}
\label{sec:setup}
 Consider a domain whose length is much greater than its width. As mentioned earlier, we only
 focus on such domains in this paper. Let us suppose that the domain is divided along its
 length into a number of non-overlapping sub-domains. Also, assume that the number of
 sub-domains is a power of two; this assumption is made mainly for ease of exposition and 
 implementation. 
 
 We organize the sub-domains into a nested hierarchy using a binary tree data structure \cite{clr90}. 
 Each node of the tree has at most two child nodes. The root of the binary tree is at level 0 of the
 tree and the children of any node in the tree are at one level greater than that node's level. The height of
 the tree is the greatest level in the tree.
  
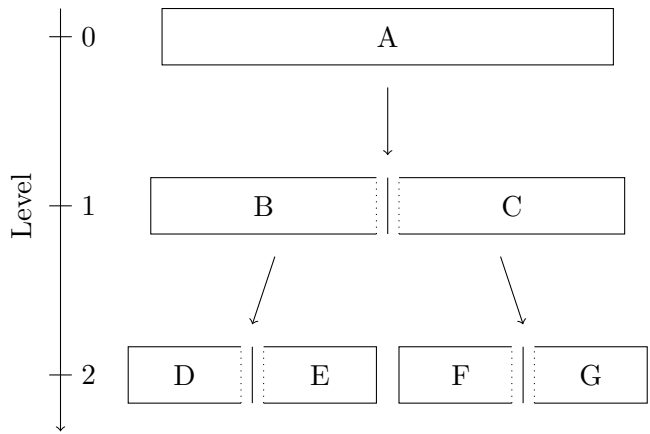
\begin{figure}
\begin{center}
\input{binaryTreeDecomposition}
\caption{\label{fig:decomposition} Binary tree based domain decomposition. A is the root of the tree that represents the full
domain. B and C are pseudo sub-domains. D, E, F and G are the true sub-domains. The height of the tree is 2.}
\end{center}
\end{figure}

 The root of the tree represents the entire domain. First, we divide the domain equally along its 
 length into two non-overlapping sub-domains. These sub-domains form the left and right children
 of the root. Each of these sub-domains is further split into two non-overlapping sub-domains in
 a similar manner; these sub-domains form the next level of the tree. The proceedure is applied 
 recursively until a pre-determined level. The sub-domains corresponding to the leaves (nodes at the last
 level) of the tree are the ``true'' sub-domains that we started out with. The sub-domains corresponding
 to the other nodes of the tree are ``pseudo'' sub-domains. Each interior (non-leaf) node of the tree
 also represents an interface between two sub-domains. Figure \ref{fig:decomposition} illustrates this
 hierarchical domain decompositon. 
  
%% file: binaryTreeDecomposition.tex
\begin{tikzpicture}[scale=1.5]
\draw (0, 0) rectangle (4, 0.5);
\draw (2, 0.25) node {A};
\draw (0.9, -1.25) node {B};
\draw (3.1, -1.25) node {C};
\draw (0.2, -2.75) node {D};
\draw (1.4, -2.75) node {E};
\draw (2.65, -2.75) node {F};
\draw (3.8, -2.75) node {G};

\draw[->] (2, -0.2) -- ++(0, -0.6);

\draw (2, -1.5) -- ++(0, 0.5);

\draw[dotted] (1.9, -1.5) -- ++(0, 0.5);
\draw (1.9, -1) -- ++(-2, 0) -- ++(0, -0.5) -- ++(2, 0);

\draw[dotted] (2.1, -1.5) -- ++(0, 0.5);
\draw (2.1, -1) -- ++(2, 0) -- ++(0, -0.5) -- ++(-2, 0);

\draw[->] (1, -1.7) -- (0.8, -2.3);

\draw (0.8, -3.0) -- ++(0, 0.5);

\draw[dotted] (0.7, -3.0) -- ++(0, 0.5);
\draw (0.7, -2.5) -- ++(-1, 0) -- ++(0, -0.5) -- ++(1, 0);

\draw[dotted] (0.9, -3.0) -- ++(0, 0.5);
\draw (0.9, -2.5) -- ++(1, 0) -- ++(0, -0.5) -- ++(-1, 0);

\draw[->] (3, -1.7) -- (3.2, -2.3);

\draw (3.2, -3.0) -- ++(0, 0.5);

\draw[dotted] (3.1, -3.0) -- ++(0, 0.5);
\draw (3.1, -2.5) -- ++(-1, 0) -- ++(0, -0.5) -- ++(1, 0);

\draw[dotted] (3.3, -3.0) -- ++(0, 0.5);
\draw (3.3, -2.5) -- ++(1, 0) -- ++(0, -0.5) -- ++(-1, 0);

\draw[->] (-0.9, 0.5) -- (-0.9, -3.25);		
\draw (-1, 0.25) -- (-0.8, 0.25);
\draw (-1, -1.25) -- (-0.8, -1.25);
\draw (-1, -2.75) -- (-0.8, -2.75);
			
\draw (-0.65, 0.25) node {$0$}
		(-0.65, -1.25) node {$1$}
			(-0.65, -2.75) node {$2$};						
		
\draw (-1.25, -1.25) node [rotate=90] {Level};				
\end{tikzpicture}

%% file: algorithm.tex
\section{Recursive Schur Decomposition Preconditioner}
\label{sec:algorithm}
Now, we will describe how we use the hierarchical domain decomposition presented in Section \ref{sec:setup} 
to solve a sparse linear system of equations represented as: $Ku = f$; $K$ is a sparse matrix, $f$ is the right hand 
side vector and $u$ is the unknown solution vector. In this paper, we will focus on systems of equations 
that are generated by applying the finite element method to discretize linear (or linearizations of nonlinear) PDEs.
We solve this problem using a Krylov subspace method (GMRES to be precise) \cite{saad2003} 
with a multi-level algorithm as preconditioner. We describe the algorithm first for the case when there are only
two sub-domains (Section \ref{sec:twoDomainCase}) and then for the general case involving several 
sub-domains (Section \ref{sec:multiDomainCase}). 

\subsection{Two Sub-domains}
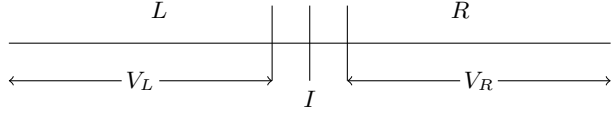
\begin{figure}
\begin{center}
\input{notations}
\caption{\label{fig:notations} Decomposition with two sub-domains.}
\end{center}
\end{figure}

Figure \ref{fig:notations} shows a decomposition of a domain into two non-overlapping sub-domains. The sub-domain on
the left is denoted by $L$ and the one on the right is denoted by $R$. The interface between the sub-domains is denoted
by $I$. The part of the sub-domain that remains after excluding the interface is denoted by $V$. The subscript on $V$
 refers to the corresponding sub-domain.

\label{sec:twoDomainCase}

Using the above notation, we can write the system of linear equations as shown in Equation (\ref{eqn:matPart}). The subscripts
$L$, $R$ and $I$ denote the degrees-of-freedom corresponding to $V_L$, $V_R$ and $I$, respectively. Note that 
$K_{II} = K_{II}^{(L)} + K_{II}^{(R)}$, where $K_{II}^{(L)}$ and $K_{II}^{(R)}$ are the contributions to $K_{II}$ from
 $L$ and $R$, respectively. 

\begin{equation}
\label{eqn:matPart}
\left [
 \begin{array}{ccc}
  K_{LL} & 0 & K_{LI} \\
   0 & K_{RR} & K_{RI} \\
   K_{IL} & K_{IR} & K_{II}
  \end{array}
 \right ]
  \left \{
 \begin{array}{l}
  u_{L} \\
   u_{R} \\
    u_{I}  
 \end{array}
  \right \}  = 
  \left \{ 
  \begin{array}{l}
   f_{L} \\ 
   f_{R} \\
    f_{I}
     \end{array} 
     \right \}
\end{equation}

Our preconditioner is based on the corresponding Schur complement system \cite{saad2003}, which is shown in Equation (\ref{eqn:schur}).

\begin{eqnarray}
\label{eqn:schur}
 S u_{I} & = & g \nonumber \\
 S  & = & K_{II} - K_{IL} {K_{LL}^{-1}} K_{LI} - K_{IR} {K_{RR}^{-1}} K_{RI} \nonumber \\
 g & = & f_{I} - K_{IL} {{K_{LL}^{-1}}} f_{L} - K_{IR} {{K_{RR}^{-1}}} f_{R}
\end{eqnarray}

Once $u_I$ is known, $u_L$ and $u_R$ can be computed as shown in Equation (\ref{eqn:backSub}).

\begin{eqnarray}
\label{eqn:backSub}
  u_{L}  & = &  {{K_{LL}^{-1}}} f_{L} - {{K_{LL}^{-1}}} K_{LI} u_{I} \nonumber \\
  u_{R} & = & {{K_{RR}^{-1}}} f_{R} - {{K_{RR}^{-1}}} K_{RI} u_{I} 
\end{eqnarray}

In our implementation, we invert the sub-domain matrices ($K_{LL}$ and $K_{RR}$) exactly using a direct method. However, 
this is not necessary and one could use an approximate solver instead. 

Note that we neither form $S$ explicitly nor solve Equation (\ref{eqn:schur}) exactly; instead, we use a few iterations
of a Krylov subspace method (GMRES to be precise) and compute an approximate solution ($\widehat{u_{I}}$). Then, we 
compute approximations to $u_{L}$ and $u_{R}$ by substituting $\widehat{u_{I}}$ for $u_{I}$ in Equation (\ref{eqn:backSub}).  
 
\subsection{Multiple Sub-domains}
\label{sec:multiDomainCase}
Next, we turn our attention to a decomposition involving several sub-domains. We represent this decomposition
using a binary tree as described in Section \ref{sec:setup}. We extend our algorithm from Section \ref{sec:twoDomainCase}
to this case by applying it to level 1 of the tree. However, we don't use a direct method to invert the matrices ($K_{LL}$ and $K_{RR}$)
 corresponding to the two pseudo sub-domains. Note that each of these pseudo sub-domains is itself a root of a sub-tree. So, we 
 could compute an approximate solution to the linear systems involving $K_{LL}$ and $K_{RR}$ by applying the algorithm 
 from Section \ref{sec:twoDomainCase} to level 1 of these sub-trees or equivalently level 2 of the original tree. This process can
 be applied recursively until we reach the last level (leaves) of the original tree. We can then invert the matrices
 ($K_{LL}$ and $K_{RR}$) corresponding to the true sub-domains exactly using a direct method or, as mentioned earlier in 
 Section \ref{sec:twoDomainCase}, using an approximate solver.

 Note that we need to invert (approximately) each of the matrices $K_{LL}$ and $K_{RR}$ at three steps in the algorithm: (1)
 computing ${{K_{LL}^{-1}}} f_{L}$ and ${{K_{RR}^{-1}}} f_{R}$ in Equation (\ref{eqn:schur}) and Equation (\ref{eqn:backSub}), (2) computing 
 ${K_{LL}^{-1}} K_{LI} v$ and ${K_{RR}^{-1}} K_{RI} v$ for some vector $v$ for evaluating the matrix-vector product (MatVec) $S v$, 
 which is required at each iteration of the Krylov solver used to solve Equation (\ref{eqn:schur}) and (3) computing
  ${K_{LL}^{-1}} K_{LI} u_{I}$ and ${K_{RR}^{-1}} K_{RI} u_{I}$ in Equation (\ref{eqn:backSub}). It will be too expensive to
  apply the overall preconditioner if we we use the recursive process described above in each of these three steps. 
  
  We would like to point out that in the second and third steps ${K_{LL}^{-1}}$ and ${K_{RR}^{-1}}$ are always applied to a
  very sparse vector because most of the rows of $K_{LI}$ and $K_{RI}$ are zero. To be precise, the rows of this vector that 
  correspond to nodes of the mesh that are not adjacent to a node on the interface are zero. We replace the recursive
  process in the second and third steps with a different approximation based on this observation. We continue to use the recursive 
  process in the first step. 
  
  Let $w$ be a sparse vector as described above, $V$ represent some pseudo sub-domain and $\widehat{V}$ represent the true 
  sub-domain contained in $V$ that shares the same interface as $V$. Consider the example shown in Figure \ref{fig:decomposition}: if
  $V$ is B then $\widehat{V}$ is E and if $V$ is C then $\widehat{V}$ is F. In the second and third steps, we approximate ${K_{VV}^{-1}} w$ 
  as follows: (1) restrict $w$ to $\widehat{V}$, (2) apply ${K_{\widehat{V}\widehat{V}}^{-1}}$ to the restricted vector and (3) prolongate
  the result back to $V$ by inserting zeros in the portion of the vector that does not belong to $\widehat{V}$.     
  
\subsection{Parallel Implementation}
\label{sec:parallel}
 In our implementation\footnote{We used the C++ programming language in our implementation.}, we assume that the
  number of true sub-domains is equal to the number of processors and we assign one true sub-domain to each
  processor. This assumption is made mainly for ease of exposition and implementation. Each interface is assigned 
  to the processor that owns the sub-domain immediately to the left\footnote{This choice is somewhat arbitrary.} of the interface.
 
 The outermost GMRES solver used for the full system, the inner GMRES solver used 
 for the Schur complement system at each level, and the direct solver used for the 
 local problems on each true sub-domain are all standard components. We used the implementation provided
 in the PETSc package \cite{petsc-home-page, petsc-user-ref, petsc-efficient} for these components. 
 
\begin{center}
\begin{pseudocode}[doublebox]{$S$-MatVec}{x}
\label{alg:smv}
\COMMENT{Computes $y = S x$} \\
\mbox{Send $x$ from $\widehat{L}$ to $\widehat{R}$} \\
\mbox{Compute $i^{(\widehat{L})} = K_{II}^{(\widehat{L})} x$} \\
\mbox{Compute $i^{(\widehat{R})} = K_{II}^{(\widehat{R})} x$} \\
\mbox{Compute $v_{\widehat{L}} = K_{\widehat{L}I} x$} \\
\mbox{Compute $v_{\widehat{R}} = K_{\widehat{R}I} x$} \\
\mbox{Solve $K_{\widehat{L}\widehat{L}} w_{\widehat{L}} = v_{\widehat{L}}$} \\
\mbox{Solve $K_{\widehat{R}\widehat{R}} w_{\widehat{R}} = v_{\widehat{R}}$} \\
\mbox{Compute $z^{(\widehat{L})} = K_{I\widehat{L}} w_{\widehat{L}}$} \\
\mbox{Compute $z^{(\widehat{R})} = K_{I\widehat{R}} w_{\widehat{R}}$} \\
\mbox{Compute $y^{(\widehat{L})} = i^{(\widehat{L})} - z^{(\widehat{L})}$} \\
\mbox{Compute $y^{(\widehat{R})} = i^{(\widehat{R})} - z^{(\widehat{R})}$} \\
\mbox{Send $y^{(\widehat{R})}$ from $\widehat{R}$ to $\widehat{L}$} \\
\mbox{Compute $y = y^{(\widehat{L})} + y^{(\widehat{R})}$} \\
\RETURN{y}
\end{pseudocode}
\end{center}
 
 The new components in the method are (1) the $S$-MatVec and (2) the Recursive Schur Decomposition (RSD)
  preconditioner. The pseudocodes for these two components are given in Algorithms \ref{alg:smv} and 
  \ref{alg:rsd}, respectively. In each of these two components, there are only two
 point-to-point communications between adjacent sub-domains; we also overlap these communications with 
 computations. We used the standard Message Passing Interface (MPI) to handle these communications.
  
\begin{center}
\begin{pseudocode}[doublebox]{RSD}{K, f, h}
\label{alg:rsd}
\COMMENT{Solves $K u = f$ approximately} \\
\COMMENT{$h$ is the height of the tree} \\
\IF h = 0
\THEN \mbox{Solve $K u = f$} 
\ELSE 
\BEGIN
\mbox{Recursion: $v_{L} \GETS$ RSD($K_{LL}$, $f_{L}$, $(h - 1)$)} \\
\mbox{Recursion: $v_{R} \GETS$ RSD($K_{RR}$, $f_{R}$, $(h - 1)$)} \\
\mbox{Compute $g^{(\widehat{L})} = K_{I\widehat{L}} v_{\widehat{L}}$} \\
\mbox{Compute $g^{(\widehat{R})} = K_{I\widehat{R}} v_{\widehat{R}}$} \\
\mbox{Send $g^{(\widehat{R})}$ from $\widehat{R}$ to $\widehat{L}$} \\
\mbox{Compute $g = f_I - g^{(\widehat{L})} - g^{(\widehat{R})}$} \\
\mbox{Solve $S u_{I} = g$ approximately} \\
\mbox{Send $\widehat{u_{I}}$ from $\widehat{L}$ to $\widehat{R}$} \\
\mbox{Compute $w_{\widehat{L}} = K_{\widehat{L}I} \widehat{u_{I}}$} \\
\mbox{Compute $w_{\widehat{R}} = K_{\widehat{R}I} \widehat{u_{I}}$} \\
\mbox{Solve $K_{\widehat{L}\widehat{L}} z_{\widehat{L}} = w_{\widehat{L}}$} \\
\mbox{Solve $K_{\widehat{R}\widehat{R}} z_{\widehat{R}} = w_{\widehat{R}}$} \\
\mbox{Compute $\widehat{u_{\widehat{L}}} = v_{\widehat{L}} - z_{\widehat{L}}$} \\
\mbox{Compute $\widehat{u_{\widehat{R}}} = v_{\widehat{R}} - z_{\widehat{R}}$} 
\END \\
\RETURN{\widehat{u}}
\end{pseudocode}
\end{center}

\subsubsection{Parallel Time Complexity}
\label{sec:complexity}
Let the number of true sub-domains be $P$, the number of unknowns per true sub-domain be $M$, and the number of 
$S$-MatVecs used to approximately solve Equation (\ref{eqn:schur}) be $\gamma$, the time (measured in terms of required 
floating point operations per processor) taken to apply the overall RSD algorithm once be $F(M, \gamma, P)$, and
 the time taken to apply the $S$-MatVec once be $G(M)$. 

$G(M)$ is dominated by the time taken to solve (probably only approximately) the true sub-domain problem once. If this is done 
using a fast algorithm such as one Multigrid V-cycle then $G(M) = \mathcal{O}(M)$ and if this is done using a direct
solver then $G(M) = \mathcal{O}(M^3)$.

Then, we have the recurrence relation: $F(M, \gamma, P) = F(M, \gamma, P/2) + \mathcal{O}(M) + \gamma G(M)$. By expanding this
 relation, we can show that $F(M, \gamma, P) = \mathcal{O}(\gamma  G(M) \log{P})$. Thus, the time taken to apply the 
 overall RSD algorithm once grows logarithmically as the number of true sub-domains.

%% file: notations.tex
\begin{tikzpicture}[scale=2]
\draw (0, 0) -- (4, 0);
\draw (2, -0.25) -- ++(0, 0.5);
\draw (1.75, -0.25) -- ++(0, 0.5);
\draw (2.25, -0.25) -- ++(0, 0.5);

\draw (2, -0.25) node[below] {\footnotesize{$I$}};
\draw (1, 0.1) node[above] {\footnotesize{$L$}};
\draw (3, 0.1) node[above] {\footnotesize{$R$}};

\draw[<-] (0, -0.25) -- (0.75, -0.25);
\draw[->] (1, -0.25) -- (1.75, -0.25);
\draw (0.875, -0.25) node {\footnotesize{$V_L$}};

\draw[<-] (2.25, -0.25) -- (3.0, -0.25);
\draw[->] (3.25, -0.25) -- (4, -0.25);
\draw (3.125, -0.25) node {\footnotesize{$V_R$}};
\end{tikzpicture}

%% file: results.tex
\section{Results and Discussion}
\label{sec:results}
 We tested our algorithm on four problems. The first is the Poisson equation shown in Equation (\ref{eqn:prob1}).
   
\begin{equation}
\label{eqn:prob1}
-\nabla \cdot \nabla u = f
\end{equation}

 The second, Equation (\ref{eqn:prob2}), and the third, Equation (\ref{eqn:prob3}), are two system PDEs where the anisotropy in
 each variable differs; these equations were previously considered in \cite{philip2012}. The coupling between 
 the $u$ and $v$ variables is weak for Equation (\ref{eqn:prob2}) and very strong for Equation (\ref{eqn:prob3}). 
 
\begin{eqnarray}
\label{eqn:prob2}
- \frac{u_{xx}}{100} - u_{yy} + \frac{v}{100} & = & f \nonumber \\
- \frac{u}{100} - v_{xx} - \frac{v_{yy}}{100} = g
\end{eqnarray}

\begin{eqnarray}
\label{eqn:prob3}
- \frac{u_{xx}}{100} - u_{yy} + 100 v & = & f \nonumber \\
- 100 u - v_{xx} - \frac{v_{yy}}{100} = g
\end{eqnarray}

 The last, Equation (\ref{eqn:prob4}), is the Navier-Lam\'{e} equation for elastostatics with 
 the first ($\lambda$) and second (shear modulus, $\mu$) Lam\'{e} parameters being 10 and 1, respectively.
 
\begin{equation}
\label{eqn:prob4}
- \nabla \cdot \nabla \vec{u} + 11 \nabla \nabla \cdot \vec{u} = \vec{f}
\end{equation}

 All the equations were solved on two dimensional rectangular domains. In all cases, homogeneous Dirichlet boundary conditions were applied on
 all sides of the domain. In this paper, we used a regular grid with $N$ nodes in each dimension to mesh each true sub-domain. The equations were 
 discretized using bi-linear finite elements. The resulting matrices were symmetric positive definite and diagonally dominant in some cases and
 unsymmetric and not diagonally dominant in others. 
 
 We performed several numerical experiements by solving each problem with various sets of parameters: $N$, $P$ and $\gamma$.
 For each problem, we picked a random solution and used the method of manufactured solutions to construct the right hand side. Starting with a zero
 initial guess, we solved each problem to a relative tolerance of 1.0e-12 in the the 2-norm of the residual. We report the number of outermost Krylov
 iterations, $\beta$, (Tables \ref{tab:prob1Iter}, \ref{tab:prob2Iter}, \ref{tab:prob3Iter} and \ref{tab:prob4Iter}) required to solve each problem as
 well as the time taken for the setup (Tables \ref{tab:prob1Setup}, \ref{tab:prob2Setup}, \ref{tab:prob3Setup} and \ref{tab:prob4Setup}) and solve 
 phases (Tables \ref{tab:prob1Solve}, \ref{tab:prob2Solve}, \ref{tab:prob3Solve} and \ref{tab:prob4Solve}) in each case. The setup phase mainly involves 
 creating the regular grid mesh, constructing the binary tree, creating the required MPI communicators, creating the matrices on each sub-domain and 
 allocating memory for the vectors. 
  
 All the experiments were performed on the {\it{Jaguar}} supercomputer at Oak Ridge National Laboratory (ORNL). The
 architectural details for this supercomputer can be found in \cite{jaguar}. 
 

\begin{table}
\footnotesize
\begin{center} 
\begin{tabular}{|c|c|c|c|c|c|c|c|c|c|c|c|c|}\hline 
 & \multicolumn{3}{c|}{P = 8} & \multicolumn{3}{c|}{P = 128} & \multicolumn{3}{c|}{P = 2048} & \multicolumn{3}{c|}{P = 32768} \\ \cline{2-13}     
 $N$ & \color{red}{$\gamma = 2$} & \color{blue}{$\gamma = 4$} & \color{green!50!black}{$\gamma = 8$} & 
   		\color{red}{$\gamma = 2$} & \color{blue}{$\gamma = 4$} & \color{green!50!black}{$\gamma = 8$} &
    	\color{red}{$\gamma = 2$} & \color{blue}{$\gamma = 4$} & \color{green!50!black}{$\gamma = 8$} &
      \color{red}{$\gamma = 2$} & \color{blue}{$\gamma = 4$} & \color{green!50!black}{$\gamma = 8$} \\ \hline 
 17 & \color{red}{18} & \color{blue}{9} & \color{green!50!black}{5} &
  		\color{red}{20} & \color{blue}{10} & \color{green!50!black}{5} & 
  		\color{red}{20} & \color{blue}{10} & \color{green!50!black}{5} &
  		\color{red}{20} & \color{blue}{10} & \color{green!50!black}{5} \\ \hline 
  33 & \color{red}{26} & \color{blue}{13} & \color{green!50!black}{6} & 
        \color{red}{27} & \color{blue}{14} & \color{green!50!black}{7} &
        \color{red}{28} & \color{blue}{14} & \color{green!50!black}{7} &
        \color{red}{27} & \color{blue}{14} & \color{green!50!black}{7} \\ \hline 
  65 & \color{red}{36} & \color{blue}{17} & \color{green!50!black}{9} &
       \color{red}{38} & \color{blue}{19} & \color{green!50!black}{9} & 
       \color{red}{38} & \color{blue}{19} & \color{green!50!black}{9} &
       \color{red}{38} & \color{blue}{19} & \color{green!50!black}{9} \\ \hline 
  129 & \color{red}{48} & \color{blue}{24} & \color{green!50!black}{12} & 
  		  \color{red}{51} & \color{blue}{25} & \color{green!50!black}{13} & 
  		  \color{red}{50} & \color{blue}{25} & \color{green!50!black}{13} &
  		  \color{red}{51} & \color{blue}{25} & \color{green!50!black}{13} \\ \hline 
 \end{tabular} 
\caption{\label{tab:prob1Iter} Number of outermost Krylov iterations required to solve Equation (\ref{eqn:prob1}).} 
\end{center} 
\end{table} 

\begin{table} 
\footnotesize
\begin{center} 
\begin{tabular}{|c|c|c|c|c|c|c|c|c|c|c|c|c|}\hline 
 & \multicolumn{3}{c|}{P = 8} & \multicolumn{3}{c|}{P = 128} & \multicolumn{3}{c|}{P = 2048} & \multicolumn{3}{c|}{P = 32768} \\ \cline{2-13}  
 $N$ & \color{red}{$\gamma = 2$} & \color{blue}{$\gamma = 4$} & \color{green!50!black}{$\gamma = 8$} & 
 			\color{red}{$\gamma = 2$} & \color{blue}{$\gamma = 4$} & \color{green!50!black}{$\gamma = 8$} &
 			\color{red}{$\gamma = 2$} & \color{blue}{$\gamma = 4$} & \color{green!50!black}{$\gamma = 8$} &
 			 \color{red}{$\gamma = 2$} & \color{blue}{$\gamma = 4$} & \color{green!50!black}{$\gamma = 8$} \\ \hline  
 17 & \color{red}{0.003} & \color{blue}{0.003} & \color{green!50!black}{0.003} &
 		  \color{red}{0.004} & \color{blue}{0.004} & \color{green!50!black}{0.004} &
 		  \color{red}{0.063} & \color{blue}{0.06} & \color{green!50!black}{0.064} &
 		  \color{red}{0.8} & \color{blue}{0.94} & \color{green!50!black}{0.881} \\ \hline 
  33 & \color{red}{0.005} & \color{blue}{0.004} & \color{green!50!black}{0.004} &
   		 \color{red}{0.006} & \color{blue}{0.006} & \color{green!50!black}{0.006} &
    	 \color{red}{0.035} & \color{blue}{0.039} & \color{green!50!black}{0.038} &
     	 \color{red}{0.934} & \color{blue}{0.808} & \color{green!50!black}{0.808} \\ \hline 
  65 & \color{red}{0.011} & \color{blue}{0.011} & \color{green!50!black}{0.011} &
   		 \color{red}{0.013} & \color{blue}{0.012} & \color{green!50!black}{0.012} &
    	 \color{red}{0.045} & \color{blue}{0.044} & \color{green!50!black}{0.044} &
       \color{red}{0.828} & \color{blue}{0.822} & \color{green!50!black}{0.837} \\ \hline 
  129 & \color{red}{0.039} & \color{blue}{0.038} & \color{green!50!black}{0.038} &
   		  \color{red}{0.04} & \color{blue}{0.04} & \color{green!50!black}{0.04} &
    		\color{red}{0.076} & \color{blue}{0.088} & \color{green!50!black}{0.076} &
     		\color{red}{0.836} & \color{blue}{1.14} & \color{green!50!black}{0.862} \\ \hline 
 \end{tabular} 
\caption{\label{tab:prob1Setup} Time (in seconds) taken for the setup phase for Equation (\ref{eqn:prob1}).} 
\end{center} 
\end{table} 
 
\begin{table} 
\footnotesize
\begin{center} 
\begin{tabular}{|c|c|c|c|c|c|c|c|c|c|c|c|c|}\hline 
 & \multicolumn{3}{c|}{P = 8} & \multicolumn{3}{c|}{P = 128} & \multicolumn{3}{c|}{P = 2048} & \multicolumn{3}{c|}{P = 32768} \\ \cline{2-13}  
 $N$ & \color{red}{$\gamma = 2$} & \color{blue}{$\gamma = 4$} & \color{green!50!black}{$\gamma = 8$} &
  \color{red}{$\gamma = 2$} & \color{blue}{$\gamma = 4$} & \color{green!50!black}{$\gamma = 8$} &
   \color{red}{$\gamma = 2$} & \color{blue}{$\gamma = 4$} & \color{green!50!black}{$\gamma = 8$} &
    \color{red}{$\gamma = 2$} & \color{blue}{$\gamma = 4$} & \color{green!50!black}{$\gamma = 8$} \\ \hline   
 17 & \color{red}{0.012} & \color{blue}{0.009} & \color{green!50!black}{0.008} & 
  \color{red}{0.02} & \color{blue}{0.01} & \color{green!50!black}{0.009} &
   \color{red}{0.018} & \color{blue}{0.012} & \color{green!50!black}{0.01} &
    \color{red}{0.033} & \color{blue}{0.07} & \color{green!50!black}{0.162} \\ \hline 
  33 & \color{red}{0.052} & \color{blue}{0.038} & \color{green!50!black}{0.029} &
   \color{red}{0.057} & \color{blue}{0.042} & \color{green!50!black}{0.034} &
    \color{red}{0.065} & \color{blue}{0.046} & \color{green!50!black}{0.037} &
     \color{red}{0.277} & \color{blue}{0.057} & \color{green!50!black}{0.044} \\ \hline 
  65 & \color{red}{0.351} & \color{blue}{0.241} & \color{green!50!black}{0.204} &
   \color{red}{0.388} & \color{blue}{0.271} & \color{green!50!black}{0.207} &
    \color{red}{0.391} & \color{blue}{0.284} & \color{green!50!black}{0.214} &
     \color{red}{0.423} & \color{blue}{0.326} & \color{green!50!black}{0.236} \\ \hline 
  129 & \color{red}{3.21} & \color{blue}{2.31} & \color{green!50!black}{1.87} &
   \color{red}{3.55} & \color{blue}{2.5} & \color{green!50!black}{2.12} &
    \color{red}{3.47} & \color{blue}{2.53} & \color{green!50!black}{2.12} &
     \color{red}{3.58} & \color{blue}{2.56} & \color{green!50!black}{2.14} \\ \hline 
 \end{tabular} 
\caption{\label{tab:prob1Solve} Time (in seconds) taken for the solve phase for Equation (\ref{eqn:prob1}).} 
\end{center} 
\end{table} 
 

\begin{table} 
\footnotesize
\begin{center} 
\begin{tabular}{|c|c|c|c|c|c|c|c|c|c|c|c|c|}\hline 
 & \multicolumn{3}{c|}{P = 8} & \multicolumn{3}{c|}{P = 128} & \multicolumn{3}{c|}{P = 2048} & \multicolumn{3}{c|}{P = 32768} \\ \cline{2-13} 
 $N$ & \color{red}{$\gamma = 2$} & \color{blue}{$\gamma = 4$} & \color{green!50!black}{$\gamma = 8$} & 
 			\color{red}{$\gamma = 2$} & \color{blue}{$\gamma = 4$} & \color{green!50!black}{$\gamma = 8$} &
  		\color{red}{$\gamma = 2$} & \color{blue}{$\gamma = 4$} & \color{green!50!black}{$\gamma = 8$} &
   		\color{red}{$\gamma = 2$} & \color{blue}{$\gamma = 4$} & \color{green!50!black}{$\gamma = 8$} \\ \hline 
 17 & \color{red}{57} & \color{blue}{25} & \color{green!50!black}{16} & 
  	  \color{red}{111} & \color{blue}{48} & \color{green!50!black}{37} &
  	  \color{red}{110} & \color{blue}{49} & \color{green!50!black}{37} &
  	  \color{red}{109} & \color{blue}{49} & \color{green!50!black}{36} \\ \hline 
  33 & \color{red}{87} & \color{blue}{36} & \color{green!50!black}{19} & 
  		 \color{red}{166} & \color{blue}{69} & \color{green!50!black}{40} &
  		  \color{red}{159} & \color{blue}{69} & \color{green!50!black}{39} &
  		   \color{red}{156} & \color{blue}{69} & \color{green!50!black}{38} \\ \hline 
  65 & \color{red}{126} & \color{blue}{52} & \color{green!50!black}{25} &
   \color{red}{249} & \color{blue}{101} & \color{green!50!black}{46} &
    \color{red}{247} & \color{blue}{100} & \color{green!50!black}{46} &
     \color{red}{245} & \color{blue}{98} & \color{green!50!black}{46} \\ \hline 
  129 & \color{red}{184} & \color{blue}{80} & \color{green!50!black}{33} &
   \color{red}{409} & \color{blue}{150} & \color{green!50!black}{63} &
    \color{red}{407} & \color{blue}{141} & \color{green!50!black}{62} &
     \color{red}{409} & \color{blue}{137} & \color{green!50!black}{60} \\ \hline 
 \end{tabular} 
\caption{\label{tab:prob2Iter} Number of outermost Krylov iterations required to solve Equation (\ref{eqn:prob2}).} 
\end{center} 
\end{table}

\begin{table} 
\footnotesize
\begin{center} 
\begin{tabular}{|c|c|c|c|c|c|c|c|c|c|c|c|c|}\hline 
 & \multicolumn{3}{c|}{P = 8} & \multicolumn{3}{c|}{P = 128} & \multicolumn{3}{c|}{P = 2048} & \multicolumn{3}{c|}{P = 32768} \\ \cline{2-13}  
 $N$ & \color{red}{$\gamma = 2$} & \color{blue}{$\gamma = 4$} & \color{green!50!black}{$\gamma = 8$} &
  		\color{red}{$\gamma = 2$} & \color{blue}{$\gamma = 4$} & \color{green!50!black}{$\gamma = 8$} &
   		\color{red}{$\gamma = 2$} & \color{blue}{$\gamma = 4$} & \color{green!50!black}{$\gamma = 8$} &
    	\color{red}{$\gamma = 2$} & \color{blue}{$\gamma = 4$} & \color{green!50!black}{$\gamma = 8$} \\ \hline 
 17 & \color{red}{0.005} & \color{blue}{0.005} & \color{green!50!black}{0.005} &
  		\color{red}{0.006} & \color{blue}{0.006} & \color{green!50!black}{0.006} &
   		\color{red}{0.048} & \color{blue}{0.038} & \color{green!50!black}{0.047} &
    	\color{red}{0.796} & \color{blue}{0.908} & \color{green!50!black}{0.96} \\ \hline 
  33 & \color{red}{0.012} & \color{blue}{0.012} & \color{green!50!black}{0.012} &
   		\color{red}{0.013} & \color{blue}{0.014} & \color{green!50!black}{0.013} &
    	\color{red}{0.045} & \color{blue}{0.045} & \color{green!50!black}{0.046} &
     	\color{red}{0.832} & \color{blue}{0.854} & \color{green!50!black}{0.805} \\ \hline 
  65 & \color{red}{0.04} & \color{blue}{0.04} & \color{green!50!black}{0.041} &
   		\color{red}{0.043} & \color{blue}{0.042} & \color{green!50!black}{0.042} &
    	\color{red}{0.149} & \color{blue}{0.096} & \color{green!50!black}{0.098} &
     	\color{red}{0.867} & \color{blue}{0.836} & \color{green!50!black}{0.862} \\ \hline 
  129 & \color{red}{0.152} & \color{blue}{0.153} & \color{green!50!black}{0.154} &
   		 \color{red}{0.157} & \color{blue}{0.158} & \color{green!50!black}{0.157} &
    	 \color{red}{0.192} & \color{blue}{0.196} & \color{green!50!black}{0.194} &
     	 \color{red}{1.09} & \color{blue}{1.05} & \color{green!50!black}{1.01} \\ \hline 
 \end{tabular} 
\caption{\label{tab:prob2Setup} Time (in seconds) taken for the setup phase for Equation (\ref{eqn:prob2}).} 
\end{center} 
\end{table} 
 
\begin{table} 
\footnotesize
\begin{center} 
\begin{tabular}{|c|c|c|c|c|c|c|c|c|c|c|c|c|}\hline 
 & \multicolumn{3}{c|}{P = 8} & \multicolumn{3}{c|}{P = 128} & \multicolumn{3}{c|}{P = 2048} & \multicolumn{3}{c|}{P = 32768} \\ \cline{2-13}  
 $N$ & \color{red}{$\gamma = 2$} & \color{blue}{$\gamma = 4$} & \color{green!50!black}{$\gamma = 8$} & 
 			 \color{red}{$\gamma = 2$} & \color{blue}{$\gamma = 4$} & \color{green!50!black}{$\gamma = 8$} &
  		 \color{red}{$\gamma = 2$} & \color{blue}{$\gamma = 4$} & \color{green!50!black}{$\gamma = 8$} &
   		 \color{red}{$\gamma = 2$} & \color{blue}{$\gamma = 4$} & \color{green!50!black}{$\gamma = 8$} \\ \hline 
 17 & \color{red}{0.051} & \color{blue}{0.033} & \color{green!50!black}{0.032} &
  		\color{red}{0.105} & \color{blue}{0.065} & \color{green!50!black}{0.078} & 
  		\color{red}{0.124} & \color{blue}{0.08} & \color{green!50!black}{0.086} &
   		\color{red}{0.203} & \color{blue}{0.356} & \color{green!50!black}{0.136} \\ \hline 
  33 & \color{red}{0.332} & \color{blue}{0.191} & \color{green!50!black}{0.159} &
   		\color{red}{0.627} & \color{blue}{0.362} & \color{green!50!black}{0.323} &
    	\color{red}{0.641} & \color{blue}{0.393} & \color{green!50!black}{0.333} &
      \color{red}{0.765} & \color{blue}{0.782} & \color{green!50!black}{0.384} \\ \hline 
  65 & \color{red}{4.05} & \color{blue}{2.42} & \color{green!50!black}{1.88} &
   		\color{red}{8.61} & \color{blue}{4.97} & \color{green!50!black}{3.63} &
    	\color{red}{8.88} & \color{blue}{4.97} & \color{green!50!black}{3.66} &
     	\color{red}{8.78} & \color{blue}{5.31} & \color{green!50!black}{3.7} \\ \hline 
  129 & \color{red}{29.8} & \color{blue}{18.7} & \color{green!50!black}{12.6} & 
  		 \color{red}{70.6} & \color{blue}{37.4} & \color{green!50!black}{25.3} &
   		 \color{red}{70.6} & \color{blue}{35.4} & \color{green!50!black}{25} &
    	 \color{red}{71.4} & \color{blue}{34.6} & \color{green!50!black}{24.4} \\ \hline 
 \end{tabular} 
\caption{\label{tab:prob2Solve} Time (in seconds) taken for the solve phase for Equation (\ref{eqn:prob2}).} 
\end{center} 
\end{table}


\begin{table} 
\footnotesize
\begin{center} 
\begin{tabular}{|c|c|c|c|c|c|c|c|c|c|c|c|c|}\hline 
 & \multicolumn{3}{c|}{P = 8} & \multicolumn{3}{c|}{P = 128} & \multicolumn{3}{c|}{P = 2048} & \multicolumn{3}{c|}{P = 32768} \\ \cline{2-13} 
 $N$ & \color{red}{$\gamma = 2$} & \color{blue}{$\gamma = 4$} & \color{green!50!black}{$\gamma = 8$} &
  		\color{red}{$\gamma = 2$} & \color{blue}{$\gamma = 4$} & \color{green!50!black}{$\gamma = 8$} &
   		\color{red}{$\gamma = 2$} & \color{blue}{$\gamma = 4$} & \color{green!50!black}{$\gamma = 8$} &
    	\color{red}{$\gamma = 2$} & \color{blue}{$\gamma = 4$} & \color{green!50!black}{$\gamma = 8$} \\ \hline 
 17 & \color{red}{32} & \color{blue}{17} & \color{green!50!black}{9} & 
 			\color{red}{34} & \color{blue}{18} & \color{green!50!black}{9} &
 			\color{red}{35} & \color{blue}{18} & \color{green!50!black}{9} &
 			\color{red}{35} & \color{blue}{18} & \color{green!50!black}{9} \\ \hline 
 33 & \color{red}{49} & \color{blue}{25} & \color{green!50!black}{12} &
  		\color{red}{51} & \color{blue}{26} & \color{green!50!black}{13} &
   		\color{red}{52} & \color{blue}{26} & \color{green!50!black}{13} &
    	\color{red}{52} & \color{blue}{26} & \color{green!50!black}{13} \\ \hline 
 65 & \color{red}{69} & \color{blue}{35} & \color{green!50!black}{17} &
  		\color{red}{74} & \color{blue}{36} & \color{green!50!black}{17} &
   		\color{red}{74} & \color{blue}{36} & \color{green!50!black}{18} &
    	\color{red}{74} & \color{blue}{36} & \color{green!50!black}{18} \\ \hline 
 129 & \color{red}{97} & \color{blue}{48} & \color{green!50!black}{23} &
  		\color{red}{102} & \color{blue}{49} & \color{green!50!black}{24} &
   	  \color{red}{104} & \color{blue}{49} & \color{green!50!black}{24} &
    	\color{red}{105} & \color{blue}{49} & \color{green!50!black}{24} \\ \hline 
 \end{tabular} 
\caption{\label{tab:prob3Iter} Number of outermost Krylov iterations required to solve Equation (\ref{eqn:prob3}).} 
\end{center} 
\end{table}

\begin{table} 
\footnotesize
\begin{center} 
\begin{tabular}{|c|c|c|c|c|c|c|c|c|c|c|c|c|}\hline 
 & \multicolumn{3}{c|}{P = 8} & \multicolumn{3}{c|}{P = 128} & \multicolumn{3}{c|}{P = 2048} & \multicolumn{3}{c|}{P = 32768} \\ \cline{2-13} 
 $N$ & \color{red}{$\gamma = 2$} & \color{blue}{$\gamma = 4$} & \color{green!50!black}{$\gamma = 8$} &
  	 	 \color{red}{$\gamma = 2$} & \color{blue}{$\gamma = 4$} & \color{green!50!black}{$\gamma = 8$} &
   		 \color{red}{$\gamma = 2$} & \color{blue}{$\gamma = 4$} & \color{green!50!black}{$\gamma = 8$} &
    	 \color{red}{$\gamma = 2$} & \color{blue}{$\gamma = 4$} & \color{green!50!black}{$\gamma = 8$} \\ \hline 
 17 & \color{red}{0.005} & \color{blue}{0.005} & \color{green!50!black}{0.005} &
  		\color{red}{0.006} & \color{blue}{0.006} & \color{green!50!black}{0.006} &
   		\color{red}{0.042} & \color{blue}{0.04} & \color{green!50!black}{0.041} &
    	\color{red}{0.832} & \color{blue}{0.848} & \color{green!50!black}{0.831} \\ \hline 
 33 & \color{red}{0.012} & \color{blue}{0.012} & \color{green!50!black}{0.012} &
  		\color{red}{0.013} & \color{blue}{0.013} & \color{green!50!black}{0.014} &
   		\color{red}{0.047} & \color{blue}{0.047} & \color{green!50!black}{0.047} &
    	\color{red}{0.801} & \color{blue}{0.861} & \color{green!50!black}{0.821} \\ \hline 
 65 & \color{red}{0.04} & \color{blue}{0.04} & \color{green!50!black}{0.04} &
 			\color{red}{0.042} & \color{blue}{0.042} & \color{green!50!black}{0.042} &
 			\color{red}{0.097} & \color{blue}{0.096} & \color{green!50!black}{0.098} &
 			\color{red}{0.815} & \color{blue}{0.843} & \color{green!50!black}{0.882} \\ \hline 
 129 & \color{red}{0.153} & \color{blue}{0.153} & \color{green!50!black}{0.153} &
 			 \color{red}{0.157} & \color{blue}{0.158} & \color{green!50!black}{0.157} &
 			 \color{red}{0.194} & \color{blue}{0.194} & \color{green!50!black}{0.193} &
 			 \color{red}{1.04} & \color{blue}{1.05} & \color{green!50!black}{1.04} \\ \hline 
 \end{tabular} 
\caption{\label{tab:prob3Setup} Time (in seconds) taken for the setup phase for Equation (\ref{eqn:prob3}).} 
\end{center} 
\end{table}

\begin{table} 
\footnotesize
\begin{center} 
\begin{tabular}{|c|c|c|c|c|c|c|c|c|c|c|c|c|}\hline 
 & \multicolumn{3}{c|}{P = 8} & \multicolumn{3}{c|}{P = 128} & \multicolumn{3}{c|}{P = 2048} & \multicolumn{3}{c|}{P = 32768} \\ \cline{2-13}  
 $N$ & \color{red}{$\gamma = 2$} & \color{blue}{$\gamma = 4$} & \color{green!50!black}{$\gamma = 8$} &
  		 \color{red}{$\gamma = 2$} & \color{blue}{$\gamma = 4$} & \color{green!50!black}{$\gamma = 8$} &
   		 \color{red}{$\gamma = 2$} & \color{blue}{$\gamma = 4$} & \color{green!50!black}{$\gamma = 8$} &
    	 \color{red}{$\gamma = 2$} & \color{blue}{$\gamma = 4$} & \color{green!50!black}{$\gamma = 8$} \\ \hline  
 17 & \color{red}{0.03} & \color{blue}{0.023} & \color{green!50!black}{0.02} &
  	  \color{red}{0.045} & \color{blue}{0.027} & \color{green!50!black}{0.021} &
   		\color{red}{0.041} & \color{blue}{0.031} & \color{green!50!black}{0.024} &
    	\color{red}{0.096} & \color{blue}{0.05} & \color{green!50!black}{0.071} \\ \hline 
 33 & \color{red}{0.187} & \color{blue}{0.139} & \color{green!50!black}{0.107} &
 			\color{red}{0.204} & \color{blue}{0.148} & \color{green!50!black}{0.121} &
 			\color{red}{0.228} & \color{blue}{0.16} & \color{green!50!black}{0.126} &
 			\color{red}{0.707} & \color{blue}{0.175} & \color{green!50!black}{0.147} \\ \hline 
 65 & \color{red}{2.28} & \color{blue}{1.67} & \color{green!50!black}{1.31} &
  		\color{red}{2.64} & \color{blue}{1.86} & \color{green!50!black}{1.42} &
   		\color{red}{2.66} & \color{blue}{1.88} & \color{green!50!black}{1.51} &
    	\color{red}{2.71} & \color{blue}{1.91} & \color{green!50!black}{1.53} \\ \hline 
 129 & \color{red}{16.2} & \color{blue}{11.6} & \color{green!50!black}{9.07} &
  		 \color{red}{18.4} & \color{blue}{12.9} & \color{green!50!black}{10.3} &
  		 \color{red}{18.8} & \color{blue}{12.9} & \color{green!50!black}{10.4} &
  		 \color{red}{19} & \color{blue}{13} & \color{green!50!black}{10.4} \\ \hline 
 \end{tabular} 
\caption{\label{tab:prob3Solve} Time (in seconds) taken for the solve phase for Equation (\ref{eqn:prob3}).} 
\end{center} 
\end{table} 
 

\begin{table} 
\footnotesize
\begin{center} 
\begin{tabular}{|c|c|c|c|c|c|c|c|c|c|c|c|c|}\hline 
 & \multicolumn{3}{c|}{P = 8} & \multicolumn{3}{c|}{P = 128} & \multicolumn{3}{c|}{P = 2048} & \multicolumn{3}{c|}{P = 32768} \\ \cline{2-13} 
 $N$ & \color{red}{$\gamma = 2$} & \color{blue}{$\gamma = 4$} & \color{green!50!black}{$\gamma = 8$} &
  		 \color{red}{$\gamma = 2$} & \color{blue}{$\gamma = 4$} & \color{green!50!black}{$\gamma = 8$} &
   		 \color{red}{$\gamma = 2$} & \color{blue}{$\gamma = 4$} & \color{green!50!black}{$\gamma = 8$} &
    	 \color{red}{$\gamma = 2$} & \color{blue}{$\gamma = 4$} & \color{green!50!black}{$\gamma = 8$} \\ \hline 
  17 & \color{red}{39} & \color{blue}{18} & \color{green!50!black}{10} &
  	   \color{red}{45} & \color{blue}{21} & \color{green!50!black}{13} &
   		 \color{red}{45} & \color{blue}{21} & \color{green!50!black}{13} &
    	 \color{red}{45} & \color{blue}{21} & \color{green!50!black}{13} \\ \hline 
  33 & \color{red}{57} & \color{blue}{26} & \color{green!50!black}{13} &
   		 \color{red}{63} & \color{blue}{29} & \color{green!50!black}{15} &
    	 \color{red}{63} & \color{blue}{29} & \color{green!50!black}{15} &
     	 \color{red}{63} & \color{blue}{29} & \color{green!50!black}{15} \\ \hline 
  65 & \color{red}{79} & \color{blue}{36} & \color{green!50!black}{17} &
   		 \color{red}{86} & \color{blue}{41} & \color{green!50!black}{19} &
    	 \color{red}{86} & \color{blue}{41} & \color{green!50!black}{19} &
     	 \color{red}{85} & \color{blue}{41} & \color{green!50!black}{19} \\ \hline 
 129 & \color{red}{112} & \color{blue}{52} & \color{green!50!black}{24} &
 			 \color{red}{119} & \color{blue}{56} & \color{green!50!black}{26} &
 			 \color{red}{118} & \color{blue}{56} & \color{green!50!black}{26} &
 			 \color{red}{118} & \color{blue}{56} & \color{green!50!black}{26} \\ \hline 
 \end{tabular} 
\caption{\label{tab:prob4Iter} Number of outermost Krylov iterations required to solve Equation (\ref{eqn:prob4}).} 
\end{center} 
\end{table}

\begin{table} 
\footnotesize
\begin{center} 
\begin{tabular}{|c|c|c|c|c|c|c|c|c|c|c|c|c|}\hline 
 & \multicolumn{3}{c|}{P = 8} & \multicolumn{3}{c|}{P = 128} & \multicolumn{3}{c|}{P = 2048} & \multicolumn{3}{c|}{P = 32768} \\ \cline{2-13} 
 $N$ & \color{red}{$\gamma = 2$} & \color{blue}{$\gamma = 4$} & \color{green!50!black}{$\gamma = 8$} &
 		   \color{red}{$\gamma = 2$} & \color{blue}{$\gamma = 4$} & \color{green!50!black}{$\gamma = 8$} &
 		   \color{red}{$\gamma = 2$} & \color{blue}{$\gamma = 4$} & \color{green!50!black}{$\gamma = 8$} &
 		   \color{red}{$\gamma = 2$} & \color{blue}{$\gamma = 4$} & \color{green!50!black}{$\gamma = 8$} \\ \hline  
 17 & \color{red}{0.005} & \color{blue}{0.005} & \color{green!50!black}{0.005} &
 		  \color{red}{0.006} & \color{blue}{0.006} & \color{green!50!black}{0.006} &
 		  \color{red}{0.041} & \color{blue}{0.067} & \color{green!50!black}{0.041} &
 		  \color{red}{0.778} & \color{blue}{0.77} & \color{green!50!black}{0.93} \\ \hline 
 33 & \color{red}{0.012} & \color{blue}{0.012} & \color{green!50!black}{0.012} &
  		\color{red}{0.013} & \color{blue}{0.013} & \color{green!50!black}{0.013} &
   		\color{red}{0.046} & \color{blue}{0.046} & \color{green!50!black}{0.046} &
    	\color{red}{0.919} & \color{blue}{0.817} & \color{green!50!black}{0.804} \\ \hline 
 65 & \color{red}{0.04} & \color{blue}{0.04} & \color{green!50!black}{0.041} &
  		\color{red}{0.042} & \color{blue}{0.043} & \color{green!50!black}{0.042} &
   		\color{red}{0.096} & \color{blue}{0.096} & \color{green!50!black}{0.098} &
    	\color{red}{0.823} & \color{blue}{0.836} & \color{green!50!black}{0.868} \\ \hline 
 129 & \color{red}{0.152} & \color{blue}{0.152} & \color{green!50!black}{0.155} &
  		 \color{red}{0.157} & \color{blue}{0.158} & \color{green!50!black}{0.158} &
   		 \color{red}{0.195} & \color{blue}{0.208} & \color{green!50!black}{0.195} &
    	 \color{red}{1.06} & \color{blue}{1.06} & \color{green!50!black}{1.05} \\ \hline 
 \end{tabular} 
\caption{\label{tab:prob4Setup} Time (in seconds) taken for the setup phase for Equation (\ref{eqn:prob4}).} 
\end{center} 
\end{table}

\begin{table} 
\footnotesize
\begin{center} 
\begin{tabular}{|c|c|c|c|c|c|c|c|c|c|c|c|c|}\hline 
 & \multicolumn{3}{c|}{P = 8} & \multicolumn{3}{c|}{P = 128} & \multicolumn{3}{c|}{P = 2048} & \multicolumn{3}{c|}{P = 32768} \\ \cline{2-13} 
 $N$ & \color{red}{$\gamma = 2$} & \color{blue}{$\gamma = 4$} & \color{green!50!black}{$\gamma = 8$} &
  		 \color{red}{$\gamma = 2$} & \color{blue}{$\gamma = 4$} & \color{green!50!black}{$\gamma = 8$} &
   		 \color{red}{$\gamma = 2$} & \color{blue}{$\gamma = 4$} & \color{green!50!black}{$\gamma = 8$} &
    	 \color{red}{$\gamma = 2$} & \color{blue}{$\gamma = 4$} & \color{green!50!black}{$\gamma = 8$} \\ \hline  
 17 & \color{red}{0.036} & \color{blue}{0.024} & \color{green!50!black}{0.022} &
 		  \color{red}{0.05} & \color{blue}{0.031} & \color{green!50!black}{0.029} &
 		  \color{red}{0.055} & \color{blue}{0.034} & \color{green!50!black}{0.033} &
 		  \color{red}{0.278} & \color{blue}{0.121} & \color{green!50!black}{0.082} \\ \hline 
 33 & \color{red}{0.217} & \color{blue}{0.144} & \color{green!50!black}{0.117} &
 		  \color{red}{0.247} & \color{blue}{0.166} & \color{green!50!black}{0.137} &
 		  \color{red}{0.297} & \color{blue}{0.175} & \color{green!50!black}{0.143} &
 		  \color{red}{0.484} & \color{blue}{0.193} & \color{green!50!black}{0.152} \\ \hline 
 65 & \color{red}{2.61} & \color{blue}{1.71} & \color{green!50!black}{1.32} &
 	 		\color{red}{3.04} & \color{blue}{2.1} & \color{green!50!black}{1.57} &
 	  	\color{red}{3.07} & \color{blue}{2.11} & \color{green!50!black}{1.59} &
 	    \color{red}{3.12} & \color{blue}{2.14} & \color{green!50!black}{1.6} \\ \hline 
 129 & \color{red}{18.4} & \color{blue}{12.4} & \color{green!50!black}{9.43} &
  		 \color{red}{21.3} & \color{blue}{14.6} & \color{green!50!black}{11.1} &
   		 \color{red}{21.2} & \color{blue}{14.7} & \color{green!50!black}{11.2} &
    	 \color{red}{21.4} & \color{blue}{14.7} & \color{green!50!black}{11.1} \\ \hline 
 \end{tabular} 
\caption{\label{tab:prob4Solve} Time (in seconds) taken for the solve phase for Equation (\ref{eqn:prob4}).} 
\end{center} 
\end{table} 

We can observe that the convergence rate of the algorithm does not deteriorate as the number of sub-domains increases. Theoretically,
 we estimated the time taken to apply the multi-level algorithm once to grow logarithmically with $P$ but, in practice the time
  taken for this seems to be roughly independent of $P$. This suggests that the constants in the complexity estimate are small.

In this work, we did not use any preconditioner for solving (approximately) the Schur complement systems. As a result, the
convergence rate for the inner Krylov subspace method used to solve (approximately) these systems deteriorates as the number of
 unknowns per sub-domain increases. Consequently, if the number of inner Krylov iterations is kept fixed then the 
quality of the overall multi-level preconditioner deteriorates as $M$ increases. Hence, for a fixed $\gamma$ the number of 
outermost Krylov iterations required to solve the problem increases as $M$ increases. 

For a fixed $M$, $\beta$ appears to be inversely proportional to $\gamma$. Since the time taken to apply the overall preconditioner once
 is directly proportional to $\gamma$, one might expect the time taken for the solve phase to be unaffected by changing $\gamma$. Nevertheless,
 increasing $\gamma$ tends to reduce the time taken for the solve phase. This can be attributed to the fact that by reducing $\beta$, we are
 reducing the number of global synchronous communications (e.g., MPI reductions to compute inner-products of vectors) involving all the processors.

%% file: conclusions.tex
\section{Conclusions}
\label{sec:conclusions}
We presented a parallel algorithm using multi-level domain decomposition for solving sparse linear systems of equations. We
tested the performance of this algorithm on several PDEs. In all cases, the performance of the algorithm did not
deteriorate as the number of sub-domains increased. To our knowledge, this is the first domain decomposition algorithm that 
has scaled to such a large number of sub-domains. Also, the running time (including the setup phase) for the solver is small. 
This is a very desirable characteristic particularly in the context of solving nonlinear and transient problems as these 
involve solving several linear systems of equations. 

We can improve the performance of the overall algorithm by using a good preconditioner for the inner Krylov subspace method. Another
useful improvement would be to use more processors per sub-domain. By allowing a variable number of processors per sub-domain, we 
could make the overall algorithm well load balanced even when the distribution of the work load across the domains is not
 uniform. Finally, we could extend the algorithm to solve three dimensional problems on more general domains by working with two
 dimensional and three dimensional domain decompositions using quadtrees and octrees instead of binary trees.

%% file: acknowledgements.tex
\section*{Acknowledgements}
This work was carried out as part of the Advanced Multi-Physics (AMP) \cite{ampTM2011} code development. The development of the AMP code was 
funded by the Nuclear Energy Advanced Modeling and Simulation (NEAMS) program of the U.S. Department of Energy Office of 
Nuclear Energy, Advanced Modeling and Simulation Office. The AMP code is open-source and available by contacting the 
author, or available with a single-user license through the Radiation Safety Information Computational Center (RSICC) at
 Oak Ridge National Laboratory as Package C793.
 
This research used resources of the Oak Ridge Leadership Computing Facility at the Oak Ridge National Laboratory, which
 is supported by the Office of Science of the U.S. Department of Energy under Contract No. DE-AC05-00OR22725.